\documentclass{amsart}
\usepackage{amsmath,amssymb}

\newtheorem{theorem}{Theorem}[section]
\newtheorem{lemma}[theorem]{Lemma}

\theoremstyle{definition}

\theoremstyle{remark}
\newtheorem{remark}[theorem]{Remark}

\numberwithin{equation}{section}

\begin{document}

\title{The $\mathcal{Q}_p$ Carleson Measure Problem}

%    Information for first author
\author{Jie Xiao}
%    Address of record for the research reported here
\address{Department of Mathematics and Statistics, Memorial University of Newfoundland, St. John's, NL A1C 5S7, Canada}
%    Current address
%\curraddr{Department of Mathematics and Statistics, Case Western Reserve University, Cleveland, Ohio 43403}
\email{jxiao@math.mun.ca}
%    \thanks will become a 1st page footnote.
\thanks{Supported in part by NSERC (Canada) and Dean of Science Start-up Funds of MUN (Canada).}

%    Information for second author
%\author{Author Two}
%\address{Mathematical Research Section, School of Mathematical Sciences, Australian National University, Canberra ACT 2601, Australia}
%\email{two@maths.univ.edu.au}
%\thanks{Support information for the second author.}

%    General info

\date{}

\dedicatory{In memory of Matts Ess\'en}


\begin{abstract}
Let $\mu$ be a nonnegative Borel measure on the open unit disk
$\mathbb{D}\subset\mathbb{C}$. This note shows how to decide that
the M\"obius invariant space $\mathcal{Q}_p$, covering
$\mathcal{BMOA}$ and $\mathcal{B}$, is boundedly (resp., compactly)
embedded in the quadratic tent-type space $T^\infty_p(\mu)$.
Interestingly, the embedding result can be used to determine the
boundedness (resp., the compactness) of the Volterra-type and
multiplication operators on $\mathcal{Q}_p$.
\end{abstract}
\maketitle

\section{Introduction}

Continuing from \cite{X1}, we answer the following question, i.e.,
the Carleson measure problem for the holomorphic
$\mathcal{Q}$-spaces (which are geometric in the sense that they are
conformally invariant):

\bigskip
\noindent{\bf Question 1.1.} \textit{Let $\mu$ be a nonnegative
Borel measure on $\mathbb D$. What geometric finite (resp.,
vanishing) property must $\mu$ have in order that $\mathcal{Q}_p$ is
boundedly (resp., compactly) embedded in
$\mathcal{T}^\infty_p(\mu)$?}
\bigskip

Here, $\mathbb D=\{z\in\mathbb C: |z|<1\}$ and $p\in (0,\infty)$ are
given, and $\mathcal{Q}_p$ stands for the space of all holomorphic
functions $f$ on $\mathbb D$ satisfying
$$
\|f\|_{\mathcal{Q}_p}=|f(0)|+\sup_{w\in\mathbb D}\sqrt{\int_{\mathbb
D}|f'(z)|^2(1-|\sigma_w(z)|^2)^pdm(z)}<\infty,
$$
where $\sigma_w(z)=(w-z)/(1-\bar{w}z)$ and $dm$ are the M\"obius map
sending $w\in\mathbb C$ to $0$ and the Lebesgue area measure on
$\mathbb C$ respectively; see \cite{X2} and \cite{X3} for an
overview of the $\mathcal{Q}_p$-theory (from 1995 to 2006) -- in
particular, $\mathcal{Q}_{p_1}\subset \mathcal{Q}_{p_2}$ when
$0<p_1<p_1\le 1$ (see Aulaskari-Xiao-Zhao \cite{AXZ});
$\mathcal{Q}_1=\mathcal{BMOA}$, John-Nirenberg's BMO space of
holomorphic functions on $\mathbb D$ (see Baernstein \cite{Ba}); and
$\mathcal{Q}_p=\mathcal{B}$, $p\in (1,\infty)$, Aulaskari-Lappan's
result in \cite{AL} (including Xiao's $\mathcal{Q}_2=\mathcal{B}$ in
\cite{X0}) regarding the Bloch space $\mathcal{B}$ of all
holomorphic functions $f$ on $\mathbb D$ with $\sup_{z\in\mathbb
D}|f'(z)|(1-|z|^2)<\infty$. Meanwhile, $\mathcal{T}^\infty_p(\mu)$
denotes the quadratic tent-type space of all $\mu$-measurable
functions $f$ on $\mathbb D$ obeying
$$
\|f\|_{\mathcal{T}^\infty_p(\mu)}=\sup_{S(I)\subseteq \mathbb
D}\sqrt{|I|^{-p}\int_{S(I)}|f|^2d\mu}<\infty;
$$
where
$$
|I|=(2\pi)^{-1}\int_I |d\xi|\quad\hbox{and}\quad
S(I)=\big\{r\xi\in\mathbb D:\ \ r\in [1-|I|, 1),\ \xi\in I\big\}
$$
are the normalized length of the subarc $I$ of the unit circle
$\mathbb T=\{z\in\mathbb C: |z|=1\}$ and the Carleson square in
$\mathbb D$ respectively. In particular, $d\mu(z)=(1-|z|^2)^pdm(z)$
and $p\in (0,1]$ lead to the square tent space on $\mathbb D$
extending the disc version of the classic one ($p=1$) on the
upper-half Euclidean space discussed in \cite{CMS} and \cite{CoVe}.

Because the norm of $f\in\mathcal{Q}_p$ is comparably dominated by
the geometric quantity (see also Aulaskari-Stegenga-Xiao
\cite{ASX}):
$$
|f(0)|+\sup_{S(I)\subseteq \mathbb
D}\sqrt{|I|^{-p}\int_{S(I)}|f'(z)|^2(1-|z|^2)^p dm(z)},
$$
our answer to Question 1.1 is naturally as follows.

\begin{theorem}\label{t1.1} Let $\mu$ be a nonnegative Borel measure
on $\mathbb D$. Then the identity operator $\mathsf{I}:\
\mathcal{Q}_p\mapsto\mathcal{T}^\infty_p(\mu)$ is bounded (resp.,
compact) if and only if
$$
\|\mu\|_{\mathcal{LCM}_p}=\sup_{S(I)\subseteq\mathbb{D}}\sqrt{\frac{\mu(S(I))}{|I|^p\big(\log\frac{2}{|I|}\big)^{-2}}}<\infty\quad
\left(\hbox{resp.,}\quad\lim_{|I|\to
0}\frac{\mu(S(I))}{|I|^p\big(\log\frac{2}{|I|}\big)^{-2}}=0\right).
$$
\end{theorem}

Based on the solution to Question 1.1, we also answer the following
problem:

\bigskip
\noindent{\bf Question 1.2.} \textit{Let $g$ be holomorphic on
$\mathbb D$. What finite (resp., vanishing) property must $g$ have
in order that $\mathsf{V}_g$ or $\mathsf{U}_g$ or $\mathsf{M}_g$ is
bounded (resp., compact) on $\mathcal{Q}_p$?}
\bigskip

Here, $\mathsf{V}_g$ and $\mathsf{U}_g$ denote the Volterra-type
operators with the holomorphic symbol $g$ on $\mathbb{D}$
respectively:
$$
{\mathsf{V}}_g f(z)=\int_0^z g'(w)f(w)dw\quad\hbox{and}\quad
{\mathsf{U}}_g f(z)=\int_0^z g(w)f'(w)dw,\ \ z\in\mathbb D.
$$
At the same time, $\mathsf{M}_g$ is the pointwise multiplication
determined by
$$
\mathsf{M}_g f(z)=f(z)g(z)=f(0)g(0)+\mathsf{V}_g f(z)+\mathsf{U}_g
f(z),\ \ z\in\mathbb{D}.
$$
Of course, in the above definition $f$ is assumed to be holomorphic
on $\mathbb{D}$. Clearly, $\mathsf{V}_g f=\mathsf{U}_f g$ and this
operator generalizes the classic Ces\'aro operator
$\mathsf{C}(f)(z)=z^{-1}\int_0^z f(w)(1-w)^{-1}dw$ -- see also
Siskakis \cite{Si} for a survey on the study of such operators.

Below is the answer to Question 1.2.

\begin{theorem}\label{t1.2} Let $g$ be holomorphic on $\mathbb D$, $d\mu_{p,g}(z)=(1-|z|^2)^p|g'(z)|^2dm(z)$ and $\|g\|_{\mathcal{H}^\infty}=\sup_{z\in\mathbb{D}}|g(z)|$. Then
\item{\rm(i)} ${\mathsf{V}}_g$ is bounded (resp., compact) on $\mathcal{Q}_p$ if and only if
$\|\mu_{p,g}\|_{\mathcal{LCM}_p}<\infty$ (resp., $\lim_{|I|\to
0}|I|^{-p}\log^2 (2/|I|)\mu_{p,g}(S(I))=0$).

\item{\rm(ii)} ${\mathsf{U}}_g$ is bounded (resp., compact) on $\mathcal{Q}_p$ if and only if
$\|g\|_{\mathcal{H}^\infty}<\infty$ (resp., $g=0$).

\item{\rm(iii)} $\mathsf{M}_g$ is bounded (resp., compact) on $\mathcal{Q}_p$ if and only if $\|\mu_{p,g}\|_{\mathcal{LCM}_p}<\infty$ and
$\|g\|_{\mathcal{H}^\infty}<\infty$ (resp., $g=0$).
\end{theorem}

The proofs of Theorems 1.1-1.2 will be given in the subsequent two
sections where the symbol $U\approx V$ will mean that there are two
constants $c_1$ and $c_2$ independent of said or implied variables
or functions such that $c_1V\le U\le c_2 V$, and $U\le c_2 V$ will
be simply written as $U\lesssim V$.

The author is grateful to K. Zhu for his helpful reply to a question
on the compact multiplication operator on the Bloch space, but also
to J. Pau and J. A. Pelaez for their valuable discussions on the
$\mathcal{Q}_p$-multiplier conjecture posed in \cite{X1}.

\section{Proof of Theorem \ref{t1.1}}

Suppose the statement before the if and only if of Theorem
\ref{t1.1} is true. Given a subarc $I$ of $\mathbb{T}$. If
$f_w(z)=\log(1-\bar{w}z)$ where $w=(1-|I|)\zeta$ and $\zeta$ is the
center of $I$, then
$$
|f_w(z)|\approx \log(2|I|^{-1}),\quad z\in S(I) $$ and
$$
|I|^{-p}\int_{S(I)}|f_w|^2d\mu\lesssim\|f_w\|^2_{\mathcal{Q}_p}\lesssim
1.
$$
Accordingly, $\|\mu\|_{\mathcal{LCM}_p}\lesssim 1$.

Conversely, let the statement after the if and only if of Theorem
\ref{t1.1} be true. To approach the desired embedding inequality, we
recall that a nonnegative Borel measure $\nu$ on $\mathbb D$ is said
to be a Carleson measure for the weighted Dirichlet space
$\mathcal{D}_p$ of all holomorphic functions $f$ obeying
$$
\|f\|_{\mathcal{D}_p}=|f(0)|+\sqrt{\int_{\mathbb
D}|f'(z)|^2(1-|z|^2)^pdm(z)}<\infty
$$
provided $\int_{\mathbb D}|f|^2d\nu\lesssim\|f\|^2_{\mathcal{D}_p}$
-- see also Stegenga \cite{Steg}. Note that $p=1$ and $p>1$ lead to
the Carleson measure for the Hardy space
$\mathcal{H}^2=\mathcal{D}_1$ and the weighted Bergman space of all
holomorphic functions in the Lebesgue space
$\mathcal{L}^2((1-|z|^2)^{p-1}dm(z))$ with respect to
$(1-|z|^2)^{p-1}dm(z)$. The following important result (written as a
lemma for our purpose) is due to Carleson \cite{Car} (for $p=1$),
Hastings \cite{Ha} (for $p=2$), Stegenga \cite {Steg} (for $p\in
[1,\infty)$), and Arcozzi-Rochberg-Sawyer \cite{ArRoSa} (for $p\in
(0,1)$).

\begin{lemma}\label{l2.1} Let $\nu$ be a nonnegative Borel measure on
$\mathbb D$.

\item{\rm(i)} If $p\in [1,\infty)$ then $\nu$ is a Carleson measure for $\mathcal{D}_p$ if and only if
$$
\|\nu\|_{\mathcal{CMD}_p}=\sup_{S(I)\subseteq \mathbb
D}\sqrt{|I|^{-p}\nu(S(I))}<\infty.
$$

\item{\rm(ii)} If $p\in (0,1)$, then $\mu$ is a Carleson measure for $\mathcal{D}_p$ if and only if
$$
\|\nu\|_{\mathcal{CMD}_p}=\sup_{w\in\mathbb
D}\sqrt{(\nu(S(w)))^{-1}\int_{\tilde{S}(w)}(\nu(S(z)\cap
S(w)))^2(1-|z|^2)^{-p-2}dm(z)}<\infty,
$$
where
$$
S(w)=\big\{z\in\mathbb{D}:\ 1-|z|\le 1-|w|,\
|\arg(w\bar{z})|\le\pi(1-|w|)\big\}
$$
and
$$
\tilde{S}(w)=\big\{z\in\mathbb{D}:\ 1-|z|\le 2(1-|w|),\
|\arg(w\bar{z})|\le\pi(1-|w|)\big\}
$$
are the Carleson and heightened Carleson boxes with vertex at
$w\in\mathbb D$ respectively.
\end{lemma}

Since $\|\mu\|_{\mathcal{LCM}_p}<\infty$ and $\lim_{|I|\to 0}\log
({2}{|I|^{-1}})=\infty$, $\mu$ is a Carleson measure for
$\mathcal{D}_p$. This fact in the case $p\in [1,\infty)$ is clear
from Lemma \ref{l2.1} (i) because of
$\|\mu\|_{\mathcal{CMD}_p}\lesssim \|\mu\|_{\mathcal{LCM}_p}$. If
$p\in (0,1)$ then this fact is due to Pau-Pelaez \cite {PP} and
follows from Lemma \ref{l2.1}, but a short proof is included below
for completeness. Fixing a point $w\in\mathbb D$, we use Fubini's
theorem to get
\begin{eqnarray*}
&&\int_{\tilde{S}(w)}{(\mu(S(z)\cap S(w)))^2}{(1-|z|^2)^{-p-2}}dm(z)\\
&&\le\int_{\tilde{S}(w)\setminus
S(w)}\frac{\mu(S(z))\mu(S(w))}{(1-|z|^2)^{p+2}}dm(z)+\int_{S(w)}\frac{(\mu(S(z)\cap
S(w)))^2}{(1-|z|^2)^{p+2}}dm(z)\\
&&\lesssim\|\mu\|_{\mathcal{LCM}_p}^2\mu(S(w))\left(\int_{\tilde{S}(w)\setminus
S(w)}\frac{dm(z)}{(1-|z|^2)^2}+\int_{S(w)}\frac{dm(z)}{(1-|z|^2)^2\log^2\frac{2}{1-|z|^2}}\right)\\
&&\lesssim\|\mu\|_{\mathcal{LCM}_p}^2\mu(S(w))\left(1+\int_0^{|w|}\int_{\arg
w-(1-|z|)}^{\arg
w+(1-|z|)}\frac{dm(z)}{(1-|z|^2)^2\log^2\frac{2}{1-|z|^2}}\right)\\
&&\lesssim\|\mu\|_{\mathcal{LCM}_p}^2\mu(S(w))\left(1+\int_0^1\frac{dt}{(1-t)\log^2\frac{2}{1-t}}\right)\\
&&\lesssim\|\mu\|_{\mathcal{LCM}_p}^2\mu(S(w)),
\end{eqnarray*}
whence $\|\mu\|_{\mathcal{CMD}_p}\lesssim
\|\mu\|_{\mathcal{LCM}_p}<\infty$ according to Lemma \ref{l2.1}
(ii).

Given any subarc $I$ of $\mathbb T$, let $w=(1-|I|)\zeta$ and
$\zeta$ be the center of $I$. Then
$$
|f(w)|\lesssim\|f\|_{\mathcal{Q}_p}\log(2|I|^{-1}),\quad
f\in\mathcal{Q}_p
$$
and
$$
(1-|w|^2)/|1-\bar{w}z|^2\approx |I|^{-1},\quad z\in S(I).
$$
Consequently, the above-verified fact that $\mu$ is a Carleson
measure for $\mathcal{D}_p$ yields

\begin{eqnarray*}
&&|I|^{-p}\int_{S(I)}|f|^2d\mu\\
&&\lesssim |I|^{-p}\Big(\int_{S(I)}|f(z)-f(w)|^2d\mu(z)+|f(w)|^2\mu(S(I))\Big)\\
&&\lesssim (1-|w|^2)^p\int_{\mathbb{D}}\Big|\frac{f(z)-f(w)}{(1-\bar{w}z)^p}\Big|^2d\mu(z)+\|f\|^2_{\mathcal{Q}_p}\|\mu\|^2_{\mathcal{LCM}_p}\\
&&\lesssim
\|\mu\|^2_{\mathcal{LCM}_p}\left(\frac{|f(0)-f(w)|^2+\int_{\mathbb{D}}\Big|\frac{d}{dz}\Big(\frac{f(z)-f(w)}{(1-\bar{w}z)^p}\Big)\Big|^2\frac{dm(z)}{(1-|z|^2)^{-p}}
}{(1-|w|^2)^{-p}}+\|f\|^2_{\mathcal{Q}_p}\right)\\
&&\lesssim\|\mu\|^2_{\mathcal{LCM}_p}\left(\|f\|^2_{\mathcal{Q}_p}\Big(1+\frac{\big(\log\frac{2}{1-|w|^2}\big)^2}{(1-|w|^2)^{-p}}\Big)+\frac{\int_{\mathbb{D}}
\Big|\frac{d}{dz}\Big(\frac{f(z)-f(w)}{(1-\bar{w}z)^p}\Big)\Big|^2\frac{dm(z)}{(1-|z|^2)^{-p}}}{(1-|w|^2)^{-p}}\right)
\\
&&\lesssim\|\mu\|^2_{\mathcal{LCM}_p}\left(\|f\|^2_{\mathcal{Q}_p}+\int_{\mathbb{D}}\Big|\frac{f(z)-f(w)}{1-\bar{w}z}\Big|^2(1-|\sigma_w(z)|^2)^pdm(z)\right)\\
&&\lesssim\|\mu\|^2_{\mathcal{LCM}_p}\|f\|^2_{\mathcal{Q}_p}.
\end{eqnarray*}
In the last inequality we have used the following estimate:
$$
\Lambda(f,w,p)=\int_{\mathbb{D}}\Big|\frac{f(z)-f(w)}{1-\bar{w}z}\Big|^2(1-|\sigma_w(z)|^2)^pdm(z)\lesssim\|f\|_{\mathcal{Q}_p}.
$$
To check this estimate we extend largely Pau-Perlaez's argument in
\cite{PP} from $p\in (0,1)$ to $p\in (0,\infty)$. Choosing
$0<\eta<p/2$ and $u=\sigma_z(v)$, we get by Zhu's \cite[Theorem 1.12
(1)]{Zhu2} that for any $z\in\mathbb D$,
\begin{eqnarray*}
&&\int_{\mathbb{D}}\frac{(1-|u|^2)^{p-\eta}}{|1-\bar{z}u|^{2+p}|1-\bar{w}u|^2}dm(u)\\
&&=\frac{(1-|z|^2)^{-\eta}}{|1-\bar{w}z|^2}\int_{\mathbb
D}\frac{(1-|v|^2)^{p-\eta}|1-\bar{z}v|^{2\eta-p}}{|1-
\bar{v}\sigma_{z}(w)|^2}dm(v)\\
&&\lesssim\frac{(1-|z|^2)^{-\eta}}{|1-\bar{w}z|^2}\int_{\mathbb{D}}\frac{(1-|v|^2)^{\eta}}{|1-
{v}\sigma_{\bar{z}}(\bar{w})|^2}dm(v)\\
&&\lesssim\frac{(1-|z|^2)^{-\eta}}{|1-\bar{w}z|^2}.
\end{eqnarray*}
The previous estimate, together with Rochberg-Wu-Zhu's formula (see
for example \cite[(2.1)]{RoWu} and \cite[p. 75: Ex. 11]{Zh1}),
Cauchy-Bunyakovsky-Schwarz's inequality, Fubini's theorem and Zhu's
\cite[Theorem 1.12 (3)]{Zhu2}, implies a series of estimates below:
\begin{eqnarray*}
&&\Lambda(f,w,p)\\
&&=\int_{\mathbb{D}}|f\circ\sigma_w(z)-f\circ\sigma_w(0)|^2|1-\bar{w}z|^{-2}(1-|z|^2)^pdm(z)\\
&&\approx\int_{\mathbb{D}}\Big|\int_{\mathbb{D}}(f\circ\sigma_w)'(u)(1-|u|^2)^{1+p}
\Big(\frac{1-(1-\bar{u}z)^{2+p}}{\bar{z}(1-\bar{u}z)^{2+p}}\Big)dm(u)\Big|^2\frac{(1-|z|^2)^p}{|1-\bar{w}z|^2}dm(z)\\
&&\lesssim\int_{\mathbb{D}}\left(\int_{\mathbb{D}}|(f\circ\sigma_w)'(u)|^2\frac{(1-|u|^2)^{2+p+\eta}}{|1-\bar{u}z|^{2+p}}dm(u)\right)
\frac{(1-|z|^2)^{p-\eta}}{|1-\bar{w}z|^2}dm(z)\\
&&\lesssim\int_{\mathbb{D}}|(f\circ\sigma_w)'(z)|^2(1-|z|^2)^{2+p+\eta}
\left(\int_{\mathbb{D}}\frac{(1-|u|^2)^{p-\eta}}{|1-\bar{z}u|^{2+p}|1-\bar{w}u|^2}dm(u)\right)dm(z)\\
&&\lesssim\int_{\mathbb{D}}|(f\circ\sigma_w)'(z)|^2(1-|z|^2)^{2+p}|1-\bar{w}z|^{-2}dm(u)\\
&&\lesssim\|f\|_{\mathcal{Q}_p}^2,
\end{eqnarray*}
as desired.

Next, we verify the compactness part. According to Lemma 2.10 in
\cite{Tj}, it suffices to show that any bounded sequence $\{f_j\}$
in $\mathcal{Q}_p$ with $f_j\to 0$ being uniform on compacta of
$\mathbb{D}$ must obey $\|f_j\|_{\mathcal{T}^\infty_p(\mu)}\to 0$.

Assume first the vanishing condition in Theorem \ref{t1.1} holds.
For $r\in (0,1)$, define the cut-off measure
$d\mu_r=1_{\{z\in\mathbb{D}:\ r<|z|\}}d\mu$, where $1_E$ denotes the
characteristic function of a set $E\subseteq\mathbb{D}$. If $r\to 1$
then
$$
\sup_{S(I)\subseteq\mathbb{D}}\frac{\mu_r(S(I))}{|I|^p\big(\log\frac{2}{|I|}\big)^{-2}}
\lesssim \sup_{|I|\le
1-r}\frac{\mu(S(I))}{|I|^p\big(\log\frac{2}{|I|}\big)^{-2}}\to 0.
$$
Suppose $\|f_j\|_{\mathcal{Q}_p}\lesssim 1$ and $f_j\to 0$ uniformly
on compacta of $\mathbb{D}$. Then the limit $\lim_{j\to
\infty}\|f_j\|_{\mathcal{T}^\infty_p(\mu)}=0$ follows from
$$
\int_{S(I)}|f_j|^2d\mu\lesssim\int_{S(I)}|f_j|^21_{\{z\in\mathbb{D}:\
|z|\le
r\}}d\mu+\|f_j\|_{\mathcal{Q}_p}^2\log^2(2|I|^{-1})\mu_r(S(I)).
$$

On the other hand, suppose $\mathsf{I}:\ \mathcal{Q}_p\mapsto
\mathcal{T}^\infty_p(\mu)$ is compact. Let $\{I_j\}$ be a sequence
of subarcs of $\mathbb{T}$ such that $|I_j|\to 0$. If $\zeta_j$ is
the center of $I_j$, $w_j=(1-|I_j|)\zeta_j$, and
$$
f_j(z)=\Big(\log\frac{1}{1-|w_j|^2}\Big)^{-1}\Big(\log\frac{1}{1-\overline{w_j}z}\Big)^2,
$$
then $\|f_j\|_{\mathcal{Q}_p}\lesssim 1$ and $f_j\to 0$ uniformly on
compacta of $\mathbb{D}$. By the compactness of $\mathsf{I}$, we
achieve that if $j\to\infty$ then
$$
0\leftarrow\|f_j\|^2_{\mathcal{T}^\infty_p}\ge
|I_j|^{-p}\int_{S(I_j)}|f_j|^2d\mu\gtrsim
\frac{\mu(S(I_j))}{|I_j|^p\big(\log\frac{2}{|I_j|})\big)^{-2}}.
$$
In other words, the desired vanishing condition is valid.

\begin{remark}\label{r2.2} Using \cite[Theorem 6]{WuZh} or \cite[Theorem 2.5.2]{X3} we can readily prove
$$
\|f\|_{\mathcal{Q}_p}\lesssim
|f(0)|+\sqrt{\sup_{w\in\mathbb{D}}\int_{\mathbb{D}}\Big|\frac{f(z)-f(w)}{1-|z|^2}\Big|^2(1-|\sigma_w(z)|^2)^pdm(z)}.
$$
which is slightly different from the conjectured-inequality:
$$
\|f\|_{\mathcal{Q}_p}\lesssim|f(0)|+\sup_{w\in\mathbb{D}}\sqrt{\Lambda(f,w,p)}.
$$
If this last estimate is true, then a new derivative-free
characterization of $\mathcal{Q}_p$ is discovered.
\end{remark}

\section{Proof of Theorem \ref{t1.2}}

(i) Note that $(\mathsf{V}_g f)'(z)=f(z)g'(z)$. So, the boundedness
part of Theorem \ref{t1.1} implies that $\mathsf{V}_g$ maps
boundedly $\mathcal{Q}_p$ into itself is equivalent to
$\|\mu_{p,g}\|_{\mathcal{LCM}_p}<\infty$, as desired. The
corresponding compactness can be demonstrated similarly.
Nevertheless, in the sequel we provide a different argument which
seems to be of independent interest. We begin with the following
density result.

\begin{lemma}\label{l3.1} If $\mathcal{LQ}_p$ and
$\mathcal{LQ}_{p,0}$ denote the spaces of all holomorphic functions
$g$ on $\mathbb{D}$ such that
$\|\mu_{p,g}\|_{\mathcal{LCM}_p}<\infty$ and $\lim_{|I|\to
0}|I|^{-p}\log^2(2/|I|)\mu_{p,g}(S(I))=0$ respectively, then
$g\in\mathcal{LQ}_{p,0}$ if and only if $g\in\mathcal{LQ}_p$ and
$\lim_{r\to 1}\|\mu_{p,g-g_r}\|_{\mathcal{LCM}_p}=0$, where
$g_r(z)=g(rz)$ is the $(0,1)\ni r$-dilation of $g$.
\end{lemma}

In fact, as in \cite[Remark 1.6]{X1} we have that
$g\in\mathcal{LQ}_p$, respectively, $g\in\mathcal{LQ}_{p,0}$, is
equivalent to
$$
\sup_{w\in\mathbb{D}}\Big(\log\frac{2}{1-|w|^2}\Big)^2\int_{\mathbb{D}}\Big(\frac{1-|w|^2}{|1-\bar{w}z|^2}\Big)^pd\mu_{p,g}(z)<\infty,
$$
respectively,
$$
g\in\mathcal{LQ}_p\quad\hbox{and}\quad\lim_{|w|\to
1}\Big(\log\frac{2}{1-|w|^2}\Big)^2\int_{\mathbb{D}}\Big(\frac{1-|w|^2}{|1-\bar{w}z|^2}\Big)^pd\mu_{p,g}(z)=0.
$$
Following the arguments for \cite[Lemma 3.5]{SiZh} and \cite[Theorem
2.1]{WiX}, we can reach the assertion described in Lemma \ref{l3.1}.

Assume now $\mathsf{V}_g$ is compact on $\mathcal{Q}_p$. For each
natural number $j$ let $I_j$ be the subarc of $\mathbb{T}$ with
center $\zeta_j$ and limit $|I_j|\to 0$. If $w_j=(1-|I_j|)\zeta_j$,
then $\{w_j\}$ has a cluster point $w_0\in\mathbb{T}$. Passing to a
subsequence we may assume that $w_j\to w_0\in\mathbb{T}$ and $f_j\to
f_0$ uniformly on compacta of $\mathbb{D}$, where
$$
f_j(z)=\log(1-\overline{w_j}z)\quad\hbox{and}\quad
f_0(z)=\log(1-\overline{w_0}z).
$$
Then $\|\mathsf{V}_g(f_j-f_0)\|_{\mathcal{Q}_p}\to 0$. As in the
proof of Theorem \ref{t1.1}, we obtain
\begin{eqnarray*}
&&|I_j|^{-p}\log^2(2/|I_j|)\mu_{p,g}(S(I_j))\\
&&\lesssim|I_j|^{-p}\int_{S(I_j)}|f_j(z)|^2 d\mu_{p,g}(z)\\
&&\lesssim |I_j|^{-p}\left(\int_{S(I_j)}|f_j(z)-f_0(z)|^2 d\mu_{p,g}(z)+\int_{S(I_j)}|f_0(z)|^2d\mu_{p,g}(z)\right)\\
&&\lesssim\|\mathsf{V}_g(f_j-f_0)\|_{\mathcal{Q}_p}^2+|I_j|^{-p}\int_{S(I_j)}|f_0(z)|^2d\mu_{p,g}(z).
\end{eqnarray*}
In the meantime, if $\mathcal{Q}_{p,0}$ stands for the space of all
functions
$$
f\in\mathcal{Q}_p\quad\hbox{with}\quad \lim_{|w|\to
1}\int_{\mathbb{D}}|f'(z)|^2(1-|\sigma_w(z)|^2)^pdm(z)=0,
$$
then $\mathcal{Q}_{p,0}$ is a closed subspace of $\mathcal{Q}_p$ and
the second dual of $\mathcal{Q}_{p,0}$ is isomorphic to
$\mathcal{Q}_p$ under the Cauchy pairing (see \cite{AlCaSi}):
$$
\langle f,h\rangle=(2\pi)^{-1}\lim_{r\to
1}\int_{\mathbb{T}}f(r\zeta)\overline{h(r\zeta)}|d\zeta|.
$$
hence $\mathsf{V}_g$ is compact and so weakly compact on
$\mathcal{Q}_{p,0}$. If $\mathsf{V}_g^\ast$ denotes the adjoint of
$\mathsf{V}_g$ acting on the dual space $\mathcal{Q}_{p,0}^\ast$ of
$\mathcal{Q}_{p,0}$ under $\langle\cdot,\cdot\rangle$, then the
adjoint $\mathsf{V}_g^{\ast\ast}$ of $\mathsf{V}_g^\ast$ is not only
bounded, but also equal $\mathsf{V}_g$ on $\mathcal{Q}_p$ -- which
can be established through the following formula for
$f\in\mathcal{Q}_{p,0}$ and $h\in\mathcal{Q}_{p,0}^\ast$:
$$
\langle \mathsf{V}_g(f),h\rangle=\langle
f,\mathsf{V}_g^\ast(h)\rangle=\overline{\langle
\mathsf{V}_g^\ast(h),f\rangle}=\overline{\langle
h,\mathsf{V}_g^{\ast\ast}(f)\rangle}=\langle\mathsf{V}_g^{\ast\ast}(f),h\rangle.
$$
Accordingly, $\mathsf{V}_g$ maps $\mathcal{Q}_p$ into
$\mathcal{Q}_{p,0}$. This implies $\mathsf{V}_g
(f_0)\in\mathcal{Q}_{p,0}$ and so
$$
\lim_{j\to\infty}|I_j|^{-p}\int_{S(I_j)}|f_0(z)|^2d\mu_{p,g}(z)=0.
$$
Consequently,
$$
\lim_{j\to\infty}|I_j|^{-p}\log^2(2/|I_j|)\mu_{p,g}(S(I_j))=0.
$$
This actually means $g\in\mathcal{LQ}_{p,0}$.

On the other hand, suppose $g\in\mathcal{LQ}_{p,0}$. Using the
uniform continuity of $g'_r$ on the closed unit disk
$\overline{\mathbb{D}}$ and Lemma \ref{l3.1} we obtain a sequence of
polynomials $\{p_k\}$ such that
$\lim_{k\to\infty}\|\mu_{p,g-p_k}\|_{\mathcal{LCM}_p}=0$. It is not
hard to see that $\mathsf{V}_{p_k}$ is compact on $\mathcal{Q}_p$.
As a matter of fact, let $\{f_j\}$ be any sequence with
$\|f_j\|_{\mathcal{Q}_p}\le 1$ and $f_j\to 0$ uniformly on compacta
of $\mathbb{D}$. Then for the polynomial $p_k$, the number $r\in
(0,1)$, and the cut-off measure
$dm_{p,r}(z)=(1-|z|^2)^p1_{\{z\in\mathbb{D}:\ |z|>r\}}dm(z)$, we use
the boundedness part of Theorem \ref{t1.1} to obtain
\begin{eqnarray*}
&&\|\mathsf{V}_{p_k}f_j\|_{\mathcal{Q}_p}^2\\
&&\lesssim\|p_k'\|_{\mathcal{H}^\infty}^2\sup_{S(I)\subseteq\mathbb{D}}|I|^{-p}\int_{S(I)}|f_j(z)|^2(1-|z|^2)^pdm(z)\\
&&\lesssim\|p_k'\|^2_{\mathcal{H}^\infty}\Big(\int_{\{z\in\mathbb{D}:\
|z|\le
r\}}|f_j(z)|^2dm(z)+\sup_{S(I)\subseteq\mathbb{D}}|I|^{-p}\int_{S(I)}|f_j(z)|^2dm_{p,r}(z)\Big)\\
&&\lesssim
\|p_k'\|^2_{\mathcal{H}^\infty}\Big(\int_{\{z\in\mathbb{D}:\ |z|\le
r\}}|f_j(z)|^2dm(z)+\|f_j\|_{\mathcal{Q}_p}^2\|m_{p,r}\|_{\mathcal{LCM}_p}^2\Big)\\
&& \lesssim
\|p_k'\|^2_{\mathcal{H}^\infty}\Big(\int_{\{z\in\mathbb{D}:\ |z|\le
r\}}|f_j(z)|^2dm(z)+(1-r)\sup_{S(I)\subseteq\mathbb{D}}|I|\log^2(2/|I|)\Big).
\end{eqnarray*}
Note that to $\epsilon>0$ there corresponds an $r_0\in (0,1)$ and a
natural number $N$ such that $j>N$ implies
$$
1-r_0<\epsilon\quad\hbox{and}\quad \int_{\{z\in\mathbb{D}:\ |z|\le
r_0\}}|f_j(z)|^2dm(z)<\epsilon.
$$
Accordingly,
$\lim_{j\to\infty}\|\mathsf{V}_{p_k}f_j\|_{\mathcal{Q}_p}=0$,
yielding the compactness of $\mathsf{V}_{p_k}:\
\mathcal{Q}_p\mapsto\mathcal{Q}_p$. Therefore, the following
estimate:
$$
\|\mathsf{V}_g-\mathsf{V}_{p_k}\|=\|\mathsf{V}_{g-p_k}\|\lesssim\|\mu_{p,g-p_k}\|_{\mathcal{LCM}_p}
$$
derives the compactness of $\mathsf{V}_g$ on $\mathcal{Q}_p$.

(ii) If $\|g\|_{\mathcal{H}^\infty}<\infty$, then the boundedness of
$\mathsf{U}_g$ follows from
$$
\|\mathsf{U}_g
f\|_{\mathcal{Q}_p}\lesssim\|g\|_{\mathcal{H}^\infty}\|f\|_{\mathcal{Q}_p},\quad
f\in\mathcal{Q}_p.
$$
Conversely, suppose $\mathsf{U}_g$ is bounded on $\mathcal{Q}_p$.
Then its operator norm $\|\mathsf{U}_g\|<\infty$. Given a nonzero
point $w\in\mathbb D$, there exists a Carleson square $S(I)$ such
that
$$
\{z\in\mathbb{D}: |\sigma_w(z)|<2^{-1}\}\subset
S(I)\quad\hbox{and}\quad 1-|w|^2\approx |I|.
$$
Choosing $f_w(z)=(\bar{w})^{-1}\log(1-\bar{w}z)$ we employ the
boundedness of $\mathsf{U}_g$ to obtain
$\|f_w\|_{\mathcal{Q}_p}\approx 1$ and
\begin{eqnarray*}
&&|I|^p|g(w)|^2\\
&&\lesssim\int_{\{z\in\mathbb{D}:
|\sigma_w(z)|<1/2\}}\frac{|g(z)|^2(1-|z|^2)^p}{|1-\bar{w}z|^2}dm(z)\\
&&\lesssim\int_{S(I)}\frac{|g(z)|^2(1-|z|^2)^p}{|1-\bar{w}z|^2}dm(z)\\
&& \lesssim\|\mathsf{U}_g(f_w)\|^2_{\mathcal{Q}_p}|I|^p\\
&&\lesssim\|\mathsf{U}_g\|^2|I|^p,
\end{eqnarray*}
and consequently, $g\in\mathcal{H}^\infty$.

As with the compactness, it is enough to verify that if
$\mathsf{U}_g:\ \mathcal{Q}_p\mapsto \mathcal{Q}_p$ is compact then
$g=0$. According to the boundedness part of (ii), the compactness of
$\mathsf{U}_g$ on $\mathcal{Q}_p$ implies $g\in\mathcal{H}^\infty$.
Now, assume $g$ is not identically equal to $0$. Then the boundary
value function $g|_{\mathbb{T}}$ cannot be identically the zero
function thanks to the maximum principle. Accordingly, there is a
positive constant $\delta$ and a sequence $\{w_j\}$ in $\mathbb{D}$
such that $w_j\to w_0\in\mathbb{T}$ and $|g(w_j)|>\delta$. Using the
classical Schwarz's lemma for $\mathcal{H}^\infty$, we readily get
$$
|g(z_1)-g(z_2)|\le
2\|g\|_{\mathcal{H}^\infty}|\sigma_{z_1}(z_2)|,\quad
z_1,z_2\in\mathbb{D}.
$$
This inequality implies that there is a sufficiently small number
$r>0$ such that $|g(z)|\ge\delta/2$ for all $j$ and $z$ obeying
$|\sigma_{w_j}(z)|<r$. Note that each pseudo-hyperbolic ball
$\{z\in\mathbb{D}:\ |\sigma_{w_j}(z)|<r\}$ is contained in a
Carleson box $S(I_j)$ with $|I_j|\approx 1-|w_j|^2$. So, if
$$
f_j(z)=\log(1-\overline{w_j}z)\quad\hbox{and}\quad
f_0(z)=\log(1-\overline{w_0}z),
$$
then $f_j-f_0\to 0$ uniformly on compacta of $\mathbb{D}$. By making
a change of variable $z=\sigma_{w_j}(u)$ and noticing
$|\sigma_{w_0}(w_j)|=1$, we obtain a series of estimates below:
\begin{eqnarray*}
&&\|\mathsf{U}_g(f_j-f_0)\|^2_{\mathcal{Q}_p}\\
&&\gtrsim |I_j|^{-p}\int_{S(I_j)}|f_j'(z)-f_0'(z)|^2|g(z)|^2(1-|z|^2)^p dm(z)\\
&&\approx |I_j|^{-p}\int_{S(I_j)}\Big|\frac{\overline{w_j}-\overline{w_0}}{(1-\overline{w_j}z)(1-\overline{w_0}z)}\Big|^2|g(z)|^2(1-|z|^2)^p dm(z)\\
&&\gtrsim \delta^2|I_j|^{-p}\int_{\{z\in\mathbb{D}:
|\sigma_{w_j}(z)|<r\}}
\Big|\frac{\overline{w_j}-\overline{w_0}}{(1-\overline{w_j}z)(1-\overline{w_0}z)}\Big|^2(1-|z|^2)^pdm(z)\\
&&\gtrsim\delta^2|I_j|^{-p}(1-|w_j|^2)^p\int_{\{u\in\mathbb{D}:\ |u|<r\}}\frac{(1-|u|^2)^p|\sigma_{w_0}(w_j)|^2}{|1-\overline{w_j}u|^{2p}|1-u\sigma_{\overline{w_j}}(\overline{w_0})|^2}dm(u)\\
&&\gtrsim\delta^2\int_{\{u\in\mathbb{D}:\ |u|<r\}}(1-|u|)^pdm(u).
\end{eqnarray*}
However, the compactness of $\mathsf{U}_g$ on $\mathcal{Q}_p$ forces
$\|\mathsf{U}_g(f_{j}-f_0)\|^2_{\mathcal{Q}_p}\to 0$, and
consequently, $\delta=0$, contradicting $\delta>0$. Therefore, $g$
must be the zero function.

(iii) The ``if" part follows immediately from the corresponding ones
of (i) and (ii). To see the ``only if" part, note that
$f_w(z)=\log(2/(1-\bar{w}z))$ belongs to $\mathcal{Q}_p$ uniformly,
i.e., $\|f_w\|_{\mathcal{Q}_p}\lesssim 1$ and any function
$f\in\mathcal{Q}_p$ has the growth
$$
|f(z)|\lesssim\|f\|_{\mathcal{Q}_p}\log\frac{2}{1-|z|^2},\quad
z\in\mathbb{D}.
$$
So, if $\mathsf{M}_g$ is bounded on $\mathcal{Q}_p$, then for fixed
$w\in \mathbb{D}$,
$$
|f_w(z)g(z)|\lesssim\|\mathsf{M}_g
f_w\|_{\mathcal{Q}_p}\log\frac{2}{1-|z|^2}\lesssim \|\mathsf{M}_g
\|\log\frac{2}{1-|z|^2},\quad z\in\mathbb{D},
$$
and hence $|g(w)|\lesssim\|\mathsf{M}_g\|$ (upon taking $z=w$ in the
last estimate), that is, $\|g\|_{\mathcal{H}^\infty}<\infty$ --
equivalently, $\mathsf{U}_g$ is bounded on $\mathcal{Q}_p$ by the
boundedness part of (ii). Consequently, $\mathsf{V}_g f=\mathsf{M}_g
f-f(0)g(0)-\mathsf{U}_g f$ gives the boundedness of $\mathsf{V}_g$
on $\mathcal{Q}_p$ and then
$\|\mu_{p,g}\|_{\mathcal{LCM}_p}<\infty$.

Suppose now $\mathsf{M}_g$ is compact on $\mathcal{Q}_p$. Then this operator is bounded and hence $\|g\|_{\mathcal{H}^\infty}<\infty$. For any nonzero sequence $\{w_j\}$ in $\mathbb{D}$ let
$$
f_j(z)=\Big(\log\frac{1}{1-|w_j|^2}\Big)^{-1}\Big(\log\frac{1}{1-\overline{w_j}z}\Big)^2.
$$
Assume $|w_j|\to 1$. Then $\|f_j\|_{\mathcal{Q}_p}\lesssim 1$ and $f_j\to 0$ uniformly on any compacta of $\mathbb{D}$. So, $\|\mathsf{M}_g(f_j)\|_{\mathcal{Q}_p}\to 0$. Because of
$$
|g(z)f_j(z)|=|\mathsf{M}_g(f_j)(z)|\lesssim\|\mathsf{M}_g(f_j)\|_{\mathcal{Q}_p}\log\frac{2}{1-|z|^2},\quad z\in\mathbb{D},
$$
we get (by letting $z=w_j$)
$$
|g(w_j)|\log\frac{1}{1-|w_j|^2}\lesssim \|\mathsf{M}_g(f_j)\|_{\mathcal{Q}_p}\log\frac{2}{1-|w_j|^2},
$$
whence $g(w_j)\to 0$. Since $g$ is bounded holomorphic function on
$\mathbb{D}$, it follows that $g=0$. We are done.

\begin{remark}\label{r3.1} The boundedness part of $\mathsf{M}_g$, plus Xiao's Theorem 1.3 (i)
and Corollary 1.4 in \cite{X1}, implies that
$$
g\in\mathsf{M}(Q_p)=\{g\in\mathcal{Q}_p:\
\mathsf{M}_g(\mathcal{Q}_p)\subseteq\mathcal{Q}_p\}
$$
if and only if
$$
g\in\mathcal{H}^\infty\ \ \hbox{and}\ \ \sup_{w\in\mathbb
D}\Big(\log\frac{2}{1-|w|^2}\Big)^2\int_{\mathbb
D}|g'(z)|^2(1-|\sigma_w(z)|^2)^pdm(z)<\infty.
$$
This equivalence proves Conjecture 1.5 in \cite{X1} (where $p\in
(0,1)$) which has been verified in Pau-Pelaez's recent work
\cite{PP}. Moreover, this equivalence and \cite[Theorem 3.3]{X1}
indicate that the Carleson's corona decomposition theorem is valid
for the multipliers of $\mathcal{Q}_p$, $p\in (0,1)$; that is, for
finitely many of holomorphic functions $g_1,...,g_n$ the following
operator:
$$
\mathsf{M}_{(g_1,...,g_n)}(f_1,...,f_n)=\sum_{j=1}^n f_jg_j
$$
maps
$\mathsf{M}(\mathcal{Q}_p)\times\cdots\times\mathsf{M}(\mathcal{Q}_p)$
onto $\mathsf{M}(\mathcal{Q}_p)$ is completely determined by the
following two conditions:
$$
\inf_{z\in\mathbb{D}}\sum_{j=1}^n|g_j(z)|>0\quad\hbox{and}\quad
(g_1,...,g_n)\in\mathsf{M}(\mathcal{Q}_p)\times\cdots\times\mathsf{M}(\mathcal{Q}_p).
$$
This result has been proved recently by Pau in \cite{P} using
\cite[Theorem 2.5.2]{X3}. Note that the case $p=1$ of the result is
due to Tolokonnikov \cite{To}, but the case $p>1$ remains open. Of
course, the corona type decompositions for all $\mathcal{Q}_p$ are
known; see also \cite{X1} (for $p\in (0,1)$), \cite{OrFa1} (for
$p=1$), and \cite{OrFa2} (for $p\in (1,\infty)$). Furthermore, the
results on boundedness and compactness of $\mathsf{V}_g$ on
$\mathcal{BMOA}$ and $\mathcal{B}$ can be seen from Siskakis-Zhao
\cite{SiZh}, Zhao \cite{Zha} and MacCluer-Zhao \cite{MaZh}; while
the boundedness result of $\mathsf{U}_g$ acting on $\mathcal{BMOA}$
is due to Danikas \cite{Da}. Finally, the boundedness descriptions
of $\mathsf{M}_g$ acting on $\mathcal{BMOA}$ and $\mathcal{B}$ can
be found in Stegenga \cite{Steg0}, Arazy \cite{Ar}, Brown-Shields \cite{BrSh} and Zhu \cite{Zh0}.
Meanwhile, $\mathsf{M}_g:\ \mathcal{B}\mapsto\mathcal{B}$ is never
compact unless $g=0$; see also Ohno-Zhao \cite{OhZh} and Zhu \cite{Zhu3}.
\end{remark}

\bibliographystyle{amsplain}

\begin{thebibliography}{10}

\bibitem{AlCaSi} A. Aleman, M. Carlsson and A. M. Simbotin, \textit{Preduals of $Q_p$ spaces and Carleson imbeddings of weighted Dirichlet spaces}. {Preprint}. (2005).

\bibitem {Ar} J. Arazy, \textit{Multipliers of Bloch functions},
{University of Haifa Mathematics Publication Series}. 52, 1982.

\bibitem{ArRoSa} N. Arcozzi, R. Rochberg and E. Sawyer, \textit{Carleson measures for analytic Besov spaces}. {Rev. Mat. Iberoamericana} \textbf{18} (2002), 443--510.

\bibitem {AL} R. Aulaskari and P. Lappan, \textit{Criteria for an analytic function to be Bloch and a harmonic or meromorphic function to be normal}.
Complex analysis and its applications (Hong Kong, 1993), 136--146,
Pitman Res. Notes Math. Ser., 305, Longman Sci. Tech., Harlow, 1994.

\bibitem {ASX} R. Aulaskari, D. Stegenga and J. Xiao, \textit{Some subclasses
of BMOA and their characterization in terms of Carleson measures}.
{Rocky Mountain J. Math.} \textbf{26} (1996), 485--506.


\bibitem {AXZ} R. Aulaskari, J. Xiao and R. Zhao, \textit{On subspaces and subsets of BMOA and UBC}. {Analysis} \textbf{15} (1995), 101--121.

\bibitem {Ba} A. Baernstein II, \textit{analytic functions of bounded mean oscillation}. Aspects of contemporary complex analysis
(Proc. NATO Adv. Study Inst., Univ. Durham, Durham, 1979), pp. 3--36, Academic Press, London-New York, 1980.

\bibitem {BrSh} L. Brown and A. L. Shields, \textit{Multipliers and
cyclic vectors in the Bloch space}. {Michigan Math. J.} \textbf{38}
(1991), 141--146.

%\bibitem {BoNi} B. Boe and A. Nicolau, \textit{Interpolation by functions in the Bloch space}. {J. Anal. Math.} \textbf{94} (2004),
%171--194.

\bibitem {Car} L. Carleson, \textit{Interpolations by bounded analytic functions
and the corona problem}. {Ann. of Math. (2)} \textbf{76}(1962),
547--559.

\bibitem {CoVe} W. S. Cohn and I. E. Verbitsky, \textit{Factoriztion
of tent spaces and Hankel operators}. {J. Funct. Anal.} \textbf{175}
(2000), 308--329.

\bibitem {CMS} R. R. Coifman, Y. Meyer and E. M. Stein, \textit{Some new function spaces and their applications to harmonic analysis}. {J. Funct. Anal.} \textbf{62}(1985), 304--335.

\bibitem {Da} N. Danikas, \textit{Untersuchungen uber analytsche funktionen
von beschr\"ankter mittlerer oszillation}. Vom Fachbereich 3
Mathematik der Technischen Universit\"at Berlin zur Verleihung des
akademischen Grades Doktor der Naturwissenschaften genehmigte
Dissertation, Berlin 1981.

\bibitem{Ha} W. W. Hastings, \textit{A Carleson measure theorem for Bergman spaces}. {Proc. Amer. Math. Soc.} \textbf{52} (1975), 237--241.

%\bibitem{Hu} Z. Hu, \textit{Extended Cesáro operators on the Bloch space in the unit ball of $C\sp n$}.
%Acta Math. Sci. Ser. B Engl. Ed. 23 (2003), no. 4, 561--566.

\bibitem{MaZh} B. MacCluer and R. Zhao, \textit{Vanishing logarithmic
Carleson measures}. {Illinois J. Math.} \textbf{46} (2002),
507--518.

%\bibitem{MaTj} S. Makhmutov and M. Tjani, \textit{Composition operators on some
%M\"obius invariant Banach spaces}. {Bull. Austral. Math. Soc.}
%\textbf{62} (2000), 1--19.


%\bibitem {NiX} A. Nicolau and J. Xiao, \textit{Bounded functions in Möbius invariant
%Dirichlet spaces}. J. Funct. Anal. 150 (1997), 383--425.


\bibitem {OhZh} S. Ohno and R. Zhao, {Weighted composition operators on the Bloch space}. {Bull. Austral. Math. Soc.} \textbf{63}(2001), 177--185.


\bibitem {OrFa1} J. M. Ortega and J. F\`abrega, {Pointwise multipliers and of corona type decomposition in BMOA}. {Ann. Inst. Fourier (Grenoble)} \textbf{46}(1996), 111--137.

\bibitem {OrFa2} J. M. Ortega and J. F\`abrega, {Pointwise multipliers and of decomposition theorems in analytic Besov spaces}. {Math. Z.} \textbf{235}(2000), 53--81.

\bibitem {P} J. Pau, {Multipliers of the $Q_s$ spaces and the corona theorem}. {Preprint}. (2007).

\bibitem {PP} J. Pau and J. A. Pelaez, {Multipliers of M\"obius
invariant $Q_s$ spaces}. {Preprint}. (2007).

\bibitem {RoWu} R. Rochberg and Z. Wu, \textit{A new characterization of Dirichlet type spaces and
applications}. {Illinois J. Math.} \textbf{37}(1993), 101-122.

%\bibitem {Sc} A. P. Schuster, \textit{Interpolation by Bloch functions}. Illinois
%J. Math. \textbf{43} (1999), 677--691.

%\bibitem{Se} K. Seip, \textit{Interpolation and sampling in spaces of analytic
%functions} University Lecture Series, 33. American Mathematical
%Society, Providence, RI, 2004.

\bibitem {Si} A. G. Siskakis, \textit{Semigroups of composition operators on spaces of analytic
functions, a review}. Studies on composition operators (Laramie, WY,
1996), 229--252, Contemp. Math., 213, Amer. Math. Soc., Providence,
RI, 1998.

\bibitem{SiZh} A. G. Siskakis and R. Zhao, \textit{A Volterra type operator on spaces of analytic functions}. Function spaces
(Edwardsville, IL, 1998), 299--311, Contemp. Math., 232, Amer. Math. Soc., Providence, RI, 1999.

\bibitem{Steg0} D. Stegenga, \textit{Bounded Toeplitz operators on $H\sp{1}$ and applications of the duality between $H\sp{1}$
and the functions of bounded mean oscillation}. {Amer. J. Math.}
\textbf{98}(1976), 573--589.

\bibitem{Steg} D. Stegenga, \textit{Multipliers of the Dirichlet space}. {Illinois J. Math.} \textbf{24} (1980), 113--139.

%\bibitem{Su} C. Sundberg, \textit{Values of BMOA functions on interpolating
%sequences}. {Michigan Math. J.} \textbf{31} (1984), 21--30.


\bibitem{Tj} M. Tjani, \textit{Compact composition operators on some
M\"obius invariant Banach spaces}. {PhD Thesis, Michigan State
University}, 1996.


\bibitem {To} V. A. Tolokonnikov, {The corona theorem in algebras of bounded analytic functions}. {Amer. Math. Soc. Trans.} \textbf{149}(1991), 61--93.


\bibitem {WiX} K. J. Wirths and J. Xiao, \textit{Recognizing $\mathcal{Q}_{p,0}$ functions per Dirichlet space structure}. {Bull. Belg. Math. Soc.} \textbf{8}(2001), 47--59.


\bibitem {WuZh} H. Wulan and K. Zhu, \textit{Derivative-free characterizations of $\mathcal{Q}_K$ spaces}. {J. Austral. Math. Soc. (Series A)}
\textbf{in press} (2007).

\bibitem {X0} J. Xiao, \textit{Carleson measure, atomic decomposition and free interpolation from Bloch space}. {Ann. Acad. Sci. Fenn. Ser. A. I. Math.} \textbf{19}(1994), 35--46.

\bibitem {X1} J. Xiao, \textit{The ${Q}_p$ corona theorem}. {Pacific J. Math.} \textbf{194}(2000), 491--509.

\bibitem {X2} J. Xiao, \textit{Holomorphic Q Classes}. {Lecture Notes in Math. 1767}. Springer-Verlag, Berlin, 2001.

\bibitem {X3} J. Xiao, \textit{Geometric ${Q}_p$ Functions}. {Frontiers in Math.} Birkh\"ause
Verlag, Bassel, 2006.

\bibitem {Zha} R. Zhao, \textit{On logarithmic Carleson measures}. {Acta Sci. Math. (Szeged)} \textbf{69} (2003), 605--618.

\bibitem {Zh0} K. Zhu, \textit{Multipliers of BMO in the Bergman metric with applications to Toeplitz operators}.
{J. Funct. Anal.} \textbf{87}(1989), 31--50.


\bibitem {Zh1} K. Zhu, \textit{Operator Theory in Function Spaces}. {Marcel
Dekker}. New York, 1990.

\bibitem {Zhu2} K. Zhu, \textit{Spaces of Holomorphic Functions in
the Unit Ball}. {Graduate Texts in Mathematics}. 226.
Springer-Verlag, New York, 2005.

\bibitem {Zhu3} K. Zhu, \textit{E-mail communication dated on March 23, 2007}.
\end{thebibliography}

\end{document}